\documentclass[10pt]{article}
\usepackage{amsmath,amssymb,amsthm,amscd}
\numberwithin{equation}{section}

\def\p{\partial}
\def\b{\bar}

\def\cC{{\cal C}}

\def\cA{{\mathcal A}}
\def\cC{{\mathcal C}}

\def\NN{{\mathbb N}}

\newtheorem{prop}{Proposition}[section]
\newtheorem{theo}[prop]{Theorem}
\newtheorem{lem}[prop]{Lemma}
\newtheorem{cor}[prop]{Corollary}
\newtheorem{rem}[prop]{Remark}

\newtheorem{q}[prop]{Question/Conjecture}
\def\begeq{\begin{equation}}
\def\endeq{\end{equation}}

\def\and{\quad{\rm and}\quad}

\let\lra=\longrightarrow

\def\mapright\#1{\,\smash{\mathop{\lra}\limits^{\#1}}\,}

\begin {document}
\bibliographystyle{plain}
\title{On K\"ahler manifolds with positive orthogonal bisectional curvature}

\date{ }
\author{X.X. Chen}
 % Enter your date or \today between curly braces
\maketitle

\tableofcontents
\newpage
\section{Introduction}
The famous Frankel conjecture asserts that any compact K\"ahler
manifold with positive bisectional curvature must be biholomorphic
to $\mathbb{C}\mathbb{P}^n.\;$  This conjecture was settled
affirmatively in early 1980s by  two groups of mathematicians
independently: Siu-Yau\cite{Siuy80} via differential geometry method
and Morri \cite{Mori79} by algebraic method.  There are many
interesting papers following this celebrated work; in particular to
understand the classification of K\"ahler manifold with non-negative
bisectional curvature, readers are referred to N. Mok's work
\cite{Mok88}  for further references. In 1982, R. Hamilton
\cite{Hamilton82} introduced the Ricci flow as a means to deform any
Riemannian metric in a canonical way to a Einstein metric. He
particularly showed that, in any 3-dimensional compact manifold, the
positive Ricci curvature is preserved by the Ricci flow. Moreover,
the Ricci flow deforms the metric more and more towards Einstein
metric. Consequently, he proved that the underlying manifold must be
diffeomorphic to $S^3\;$ or a finite quotient of  $S^3.\;$ By a theorem of M. Berger in the 1960s,
any K\"ahler Einstein metric with positive bisectional curvature is
the Fubni-Study metric (with constant bisectional curvature).  A
natural and long standing problem for K\"ahler Ricci flow is: In
$\mathbb{C}\mathbb{P}^n$, is K\"ahler Ricci flow converges to the
Fubni-Study metric  if the initial
metric has positive bisectional curvature? There are many interesting
work in this direction in 1990s (c.f. \cite{Bando87}  \cite{Mok88} ) and the problem was completely
settled in 2000 by \cite{chentian001} \cite{chentian002}
affirmatively. One key idea is the introduction of  a series new
geometrical functionals which play a crucial role in deriving the
bound of scalar curvature, diameter and
Sobleve constants etc..\\

  It is well known that the positivity of bisectional curvature is
  preserved along the K\"ahler Ricci flow, due to S. Bando
  \cite{Bando87} in dimension 3 and N. Mok in general dimension
  \cite{Mok88}.  Following the work of N. Mok, in an unpublished work of Cao-Hamilton, they claimed that the   positive orthogonal
  bisectional curvature is also preserved under the K\"ahler Ricci flow.\\

  In any K\"ahler manifold,  we can split the space of $(1,1)$ forms into two orthogonal components:
  the line
  spanned
  by the K\"ahler form, and its orthogonal complement $\Lambda_0^{1,1}.\;$ The traceless part of the bisectional curvature can be viewed as an operator
  acting on this subspace $\Lambda_0^{1,1}.\;$ We call a K\"ahler metric has 2-positive traceless bisectional curvature
  if the sum of any two eigenvalues
  of the tracless bisectional curvature operator in $\Lambda_0^{1,1}$ is positive.  In complex dimension 2, this property was
   studied  by \cite{[JP041]}
  where they proved that the 2-positive traceless bisectional curvature condition is preserved under the K\"ahler Ricci
  flow.  This result, as well as the method of proof, is similar to that of  H. Chen in Riemannian case.
  In \cite{chenlihao}, H. Li and the author prove that the 2-positive traceless bisectional curvature operator is preserved under the K\"ahler Ricci
  flow in all dimensions.  Moreover, we show that in \cite{chenlihao} that K\"ahler metric with 2-positive traceless bisectional curvature
  must also have positive  orthogonal bisectional curvature.
  Therefore, if the initial K\"ahler metric has 2-positive traceless bisectional curvature, then
  the orthogonal bisectional curvature remains positive on the
  entire K\"ahler Ricci flow.\\

    In this paper, we  study any K\"ahler manifold where the
positive orthogonal bisectional curvature is preserved on the
K\"ahler Ricci flow.  Naturally, we always assume that the first
Chern class $C_1$ is positive.   Under this assumption, we first
show that various inequalities (convex cone), on the scalar
curvature, the Ricci curvature tensor and the holomorphic sectional
curvature tensor, are preserved over the K\"ahler Ricci flow
respectively.  We then discuss geometrical application of these
results. In particular, we prove that any irreducible K\"ahler
manifold with positive orthogonal bisectional curvature must be
biholomorphic to $\mathbb{C}\mathbb{P}^n.\;$ This can be viewed as a
generalization of Siu-Yau\cite{Siuy80}, Morri's solution
\cite{Mori79} of the Frankel conjecture. \\

Now we state some results on maximum principle first.
\begin{theo} \label{th:scalarbecomepositive} Along the K\"ahler Ricci flow,  the following statement hold
\begin{enumerate}
\item  The lower bound of the scalar curvature, if $\leq 0$, is
preserved and improved over time.   The lower bound will increase
to $0$ exponentially. \item The upper bound of the scalar
curvature grows at most exponentially.
\end{enumerate}
\end{theo}

\begin{theo} \label{th:riccibecomepositive}
If the orthogonal bisectional curvature is positive along the
K\"ahler Ricci flow, then the following hold:
\begin{enumerate}\item If the initial metric has positive Ricci curvature, then
this condition will be preserved. \item If the initial Ricci
curvature is not positive, then the lower bound of the Ricci
curvature will increase. As $t\rightarrow \infty$, this lower
bound will at least increase to $0$ as $t\rightarrow
\infty.\;$\end{enumerate}
\end{theo}

\begin{rem} The second part of Theorem
\ref{th:scalarbecomepositive}.1 was proved firs by B. Chow
\cite{chow91} in $S^2.\;$    Statements in Theorem \ref{th:scalarbecomepositive} should be known to
the experts in the field and the proof in high dimension is similar to 2-d case.\\

Theorem \ref{th:riccibecomepositive}.1 was observed by
Cao-Hamilton first.
\end{rem}

\begin{theo} \label{th:holosectimprov} If the positivity of the orthogonal
bisectional curvature is preserved along the K\"ahler Ricci flow,
then the lower bound of the holomorphic sectional curvature, if
non-positive, is preserved under the flow.  Moreover, the lower
bound will be improved to $0$ exponentially over the flow.
\end{theo}

The following theorem is more or less technical.

\begin{theo}\label{th:positivebisectionalcurvimprove} On the K\"ahler manifold where the positivity of the
orthogonal bisectional curvature is preserved under the K\"ahler
Ricci flow, suppose that the condition $Ric\geq \nu > 0 $ is
preserved over the entire flow. Then the lower bound of the
holomorphic sectional curvature, if $\leq 0,$ will become positive
after finite time. Moreover, if the lower bound of bisectional
curvature is positive and $\nu > {1\over 2}$, then this lower bound
will be increased over time and approach ${{2\nu-1}\over
  {n+1}}$ exponentially.
\end{theo}

The following question is then very natural
\begin{q} Is an irreducible compact K\"ahler manifold with  positive orthogonal bisectional curvature
 necessary $\mathbb{C}\mathbb{P}^n$?
\end{q}

A related question is: for higher dimension, does manifold with positive orthogonal bisectional
curvature necessary has first Chern class positive? Apparently, we need to assume $n > 1.\;$
Hopefully,  when $n$ is large enough, the positive orthogonal bisectional curvature
condition implies the positivity of first Chern class?\\

Following \cite{chentian002}, we have
\begin{theo}
On a K\"ahler Einstein manifold where the positive orthogonal
bisectional curvature is preserved over the K\"ahler Ricci flow,
the flow converges to a K\"ahler Einstein metric exponentially
fast.
\end{theo}
Comparing  to Theorem 1.1 and Remark 1.9 of \cite{chentian002}, we
drop the assumption that the initial metric has positive Ricci
curvature. Theorem \ref{th:riccibecomepositive} implies that, after
finite time, the Ricci curvature will be bigger than $-1.\;$ The
energy functional $E_1$ will be monotonely decreased afterwards. The
rest of the proof is the same as in \cite{chentian002}.  Here $E_1$
is one of a set of new geometrical functional introduced in
\cite{chentian002}.

 Invoking a deep theorem of
Perelman where he shows that the scalar curvature is uniformly
bounded along the K\"ahler Ricci flow with any smooth initial
metric\footnote{Using Perelman's local estimates together with Cao's
Harnack inequalities under this case,  Cao-Zhu-Zhu \cite{[3]} gave
another proof of a uniform bound of the scalar curvature under the
assumption of positive bisectional curvature.}, then we prove

\begin{theo}\label{th:mainwithP} For any irreducible K\"ahler manifold which admits positive orthogonal bisectional
curvature and $C_1 > 0$, if this positivity condition is preserved
under the flow, then the underlying manifold is biholomorphic to
$\mathbb{C}\mathbb{P}^n.\;$
\end{theo}
As a corollary,
we have
\begin{theo}
 Any irreducible K\"ahler manifold with positive orthogonal
bisectional curvature or positive 2-traceless bisectional curvature
and $C_1 > 0$ must be $\mathbb{C}\mathbb{P}^n.\;$
\end{theo}

Following proof of Theorem 1.8 and  classification of manifold with non-negative bisectional curvature
(c.f. \cite{Mok88}), one shall be able to generalize Theorem 1.8 to the case of non-negative orthogonal bisectional curvature positive
case.\\

     One unsatisfactory feature of the present proof is that  we had to take limit as $t\rightarrow \infty$
    first in order to show that the bisectional curvature becomes positive
    after finite time.   The argument is a little bit in-direct.  If a more direct argument can be obtained,
    then perhaps one can avoid using of Perelman's  theorem on the scalar curvature function.
    The following is a weak, but direct theorem.
\begin{theo} \label{th:mainwithP2} For any K\"ahler manifold with positive orthogonal
bisectional curvature and $C_1 > 0$ ,  then the following statements
are equivalent
\begin{enumerate}
\item there exists a lower bound of the energy functional $E_1.\;$
\item there exists a lower bound of the Mabuchi energy;
 \item there is a K\"ahler Einstein metric in the K\"ahler class;
\end{enumerate}
In particular, the flow converges exponentially fast to the
Fubni-Study metric and the underlying manifold is
$\mathbb{C}\mathbb{P}^n.\;$
 \end{theo}
Note that in the proof of this theorem, we don't use Frankel
Conjecture.\\

% we have immediately

%\begin{theo} Let $(M,\omega)$ be K\"ahler manifold with positive first Chern class. Suppose
%there is a K\"ahler metric in $[\omega]$ which has non-negative orthogonal bisectional curvature.
%Then there exists non-negative integers $k, N_1, N_2,\cdots N_l$ and irreducible compact Hermitian
%symmetric spaces $M_1, M_2, \cdots M_k$ or rank $\geq 2$, such that  the universal cover of $M$  is biholomorphic to
%\[\mathbb{C}^k \times \mathbb{C} \mathbb{P}^{N_1}\times \cdots \mathbb{C}\mathbb{P}^{N_l}
%\times M_1 \cdots M_p.
%\]\end{theo}

 {\bf Acknowledgment}: The author wishes to thank S. K. Donaldson,  G. Tian
for their continuous support. The author also wishes
to thank H. Li for many interesting discussions on K\"ahler Ricci
flow and thanks for both H. Li and W. Y. He for reading the first
draft of this paper.
\section{The Maximum principle along the K\"ahler Ricci flow}
In this section, we prove some theorems on the scalar curvature,
Ricci curvature and holomorphic sectional curvature via Hamilton's
maximum principle on tensors along a geometric flow.
\subsection{On the scalar curvature}
Here we give a proof to Theorem \ref{th:scalarbecomepositive}.
\begin{proof}
 The
evolution equation for K\"ahler potential is:
\[
  {{\p \varphi}\over {\p t}} = \log {\omega_\varphi^n \over
  \omega^n} + \varphi.
\]
Set\[ F(t) = |\nabla {{\p \varphi}\over {\p t}}
\mid_{\varphi(t)}^2.
\]
Then the evolution equation for $F(t)$ is simple
\[
  {\p \over {\p t}} F = \triangle_{\varphi(t)} F + F - \left({{\p \varphi}\over {\p t}}\right)_{\alpha \b \beta} \cdot
  \left({{\p \varphi}\over {\p t}}\right)_{\b \alpha \beta} -  \left({{\p \varphi}\over {\p t}}\right)_{\alpha  \beta} \cdot
  \left({{\p \varphi}\over {\p t}}\right)_{\b \alpha \b \beta}.
\]
Note that
\[
\left({{\p \varphi}\over {\p t}}\right)_{\alpha \b \beta} =
g_{\alpha \b \beta} - R_{\alpha \b \beta}.
\]
Recalled  the evolution equation for the scalar curvature \[
\begin{array} {lcl} {{\p R}\over {\p t}} & = & \triangle R + |Ric-1|^2 + 2 R - n -
R\\ & = &   \triangle R +  \left({{\p \varphi}\over {\p
t}}\right)_{\alpha \b \beta} \cdot
  \left({{\p \varphi}\over {\p t}}\right)_{\b \alpha \beta}  + R - n.\end{array}
\]
Set $h = R-n + F.\;$ Then the  evolution equation for $h$ is
\[
\begin{array}{lcl} {{\p h}\over {\p t}} & = & \triangle h  + h - \left({{\p \varphi}\over {\p t}}\right)_{\alpha  \beta} \cdot
  \left({{\p \varphi}\over {\p t}}\right)_{\b \alpha \b \beta}.
  \end{array}
\]
Thus, the upper bound of $h$ grows at most exponentially
\[
  h \leq \displaystyle \max_{x\in M}\; h(x,0) e^{
t}.\] Consequently, we have
\[
  R(x,t) \leq C_1 e^t + n, \qquad \forall t > 0.
\]
To obtain lower bound, let $-\mu(t)$ be the negative lower bound
of the scalar curvature at time $t.\;$ Then
\[
\begin{array}{lcl} {\p\over {\p t}} (R + \mu(t))& =& \triangle (R +
\mu(t)) + Ric^2 - (R + \mu(t)) +\mu(t) + \mu(t)'\\
&\geq & \triangle (R + \mu(t)) - (R + \mu(t))
\end{array}
\]
where we set
\[
 \mu'(t) + \mu(t) = 0, \qquad \mu(0) >  - \displaystyle \min_{x\in
 M}
 R(x,0) > 0.
\]
In particular, if $\displaystyle \min_{x \in M}\; R(x,0) < 0,\;$
then
\[
  R(x,t) \geq \displaystyle \min_{x\in M} R(x,0)\; e^{-t}.
\]
\end{proof}
\subsection{On the Ricci curvature}
Now we give a proof of Theorem \ref{th:riccibecomepositive}.

\begin{proof} Suppose $-\mu(t) < 0$ is lower bound of Ricci curvature (not necessary optimal).
   Set
 \[
   \hat{R}_{i\b j} = R_{i\b j} + \mu(t) g_{i\b j}.
 \]
 Then
 \[\begin{array}{lcl} {\p \over {\p t}} \hat{R}_{i\b j} & = & {\p \over {\p t}} {R}_{i\b
 j} + \mu(t) {\p \over {\p t}} g_{i\b j} +  \mu' g_{i\b j}\\
 &= & \triangle_\varphi R_{i\b j} + R_{i \b j k\b l} R_{\b k  l} -
 R_{i \b p} R_{p \b j} + \mu(t) (g_{i\b j} - R_{i \b j}) + \mu'
 g_{i\b j}\\
 & = & \triangle_\varphi \hat{R}_{i\b j} + R_{i \b j k\b l} \hat{R}_{\b k
 l} -\mu R_{i\b j} -\hat{R}_{i\b p} \hat{R}_{p\b j} + 2 \mu
 \hat{R}_{i\b j} - \mu^2 g_{i\b j} + \mu (\mu+1) g_{i\b j} - \mu
 \hat{R}_{i\b j} +\mu' g_{i\b j}\\
& = & \triangle_\varphi \hat{R}_{i\b j} + R_{i \b j k\b l}
\hat{R}_{\b k
 l} -\mu \hat{R}_{i\b j} -\hat{R}_{i\b p} \hat{R}_{p\b j} + 2 \mu
 \hat{R}_{i\b j}  + \mu (\mu+1) g_{i\b j} - \mu
 \hat{R}_{i\b j} +\mu' g_{i\b j}\\
 & = &\triangle_\varphi \hat{R}_{i\b j} + R_{i \b j k\b l} \hat{R}_{\b k
 l}  -\hat{R}_{i\b p} \hat{R}_{p\b j} + (\mu (\mu+1) + \mu') g_{i\b
 j}.
\end{array}
 \]
If we choose $\mu(0)$ such that $R_{i\b j} + \mu g_{i\b j} \geq 0$
at time $t =0$  set $\mu(t) $ solves the following equation
\[
  \mu'(t) + \mu (\mu+1) = 0
\]
or
\[
  \mu(t) = {{C }\over { e^t- C}}, \qquad {\rm where}\; C = {\mu(0)\over {\mu(0)+1}}.
\]
Then the evolution equation for $\hat{R}_{i\b j}$ is:
\[ {\p \over {\p t}} \hat{R}_{i\b j} = \triangle_\varphi \hat{R}_{i\b j} + R_{i \b j k\b l} \hat{R}_{k \b
 l}  -\hat{R}_{i\b p} \hat{R}_{p\b j}.\]
Consequently, the non-negativity of $\hat{R}_{i\b j}$ is
preserved.  This is because at the point where where $\hat{R}_{i\b
j}$ vanishes at least in one direction, we  can show that  the
right hand side must be non-negative. In fact, set this direction
as ${\p \over {\p z_1}}$ and diagonalize the Ricci curvature at
this point. Then
\[
{\p\over {\p t}} \hat{R}_{1\b 1} \geq R_{1\b 1 j \b j} \hat{R}_{j
\b j} - \hat{R}_{1\b 1} \hat{R}_{\b 1 1} = \displaystyle
\sum_{j=2}^n \; R_{1\b 1 j \b j} \hat{R}_{j \b j} \geq 0.
\]
\end{proof}
\subsection{On the holomorphic sectional curvature}
Suppose $(M, g(0))$ is a K\"ahler manifold with positive
orthogonal bisectional curvature.  Roughly speaking, we want to
prove that the lower bound of the holomorphic sectional curvature,
if not positive, will be preserved or improved under appropriate
other conditions.  Now we give a proof to Theorem
\ref{th:holosectimprov}.

%\section{main theorem}
\begin{proof}  For any bisectional curvature type tensor
$A_{i \b j k \b l}$, we define the operator $\square$ as
\[
\begin{array}{lcl}
\square A_{i \b j k \b l} & = & \bigtriangleup A_{i \overline{j} k
\overline{l}} + A_{i \overline{j} p \overline{q}} A_{q
\overline{p} k \overline{l}} - A_{i \overline{p} k \overline{q}}
A_{p \overline{j} q \overline{l}} + A_{i \overline{l} p
\overline{q}} A_{q \overline{p} k \overline{j}} \\ & &  + A_{i
\overline{j} k \overline{l}}
  -{1\over 2} \left( R_{i \overline{p}}A_{p \overline{j} k \overline{l}}  +
R_{p \overline{j}}A_{i \overline{p} k \overline{l}} + R_{k
\overline{p}}A_{i \overline{j} p \overline{l}} + R_{p
\overline{l}}A_{i \overline{j} k \overline{p}} \right)
\end{array}
\]
Then the evolution equation for bisectional curvature is
\[
\begin{array}{lcl} {{\partial }\over {\partial t}} R_{i \overline{j} k
\overline{l}} & = & \square R_{i\b j k \b l} \\ & = &
\bigtriangleup R_{i \overline{j} k \overline{l}} + R_{i
\overline{j} p \overline{q}} R_{q \overline{p} k \overline{l}} -
R_{i \overline{p} k \overline{q}} R_{p \overline{j} q
\overline{l}} + R_{i
\overline{l} p \overline{q}} R_{q \overline{p} k \overline{j}} \\
& &  + R_{i \overline{j} k \overline{l}}
  -{1\over 2} \left( R_{i \overline{p}}R_{p \overline{j} k \overline{l}}  +
R_{p \overline{j}}R_{i \overline{p} k \overline{l}} + R_{k
\overline{p}}R_{i \overline{j} p \overline{l}} + R_{p
\overline{l}}R_{i \overline{j} k \overline{p}} \right).
\end{array}
\]
Set
\[
S_{i\b j k\b l} = R_{i \b j k \b l} - \mu(t) (g_{i\b j} g_{k \b l}
+ g_{i \b l} g_{j \b k}) = R_{i \b j k \b l} - \mu (g*g)_{i\b j k
\b l}
\]
and
\[
S_{k\b l} = R_{k \b l} - \mu (n+1) g_{k \b l}.
\]
At the point where $g_{i\b j} = \delta_{i \b j},\;$
 we can re-write the evolution equation for $R$ as
\[
\begin{array}{lcl} {{\partial }\over {\partial t}} R_{i \overline{j} k \overline{l}}
& = & \bigtriangleup \left(S_{i \overline{j} k \overline{l}} + \mu
(g*g)_{i\b j k \b l} \right) + (S_{i \overline{j} p \overline{q}}
+ \mu (g*g)_{i \b j p \b q})  (S_{q \overline{p} k \overline{l}} +
\mu (g*g)_{q\b p k \b l}) \\
& & - (S_{i \overline{p} k \overline{q}}+  \mu (g*g)_{i\b p k \b
q}) (S_{p \overline{j} q \overline{l}}+ \mu(g*g)_{p \b j q \b l})
+  S_{i \overline{j} k \overline{l}} +  \mu (g*g)_{i \b j k \b l}
\\ && + (S_{i
\overline{l} p \overline{q}}+ \mu (g*g)_{i \b l p \b q}) (S_{q \overline{p} k \overline{j}}+ \mu (g*g)_{q\b p k \b j}) \\
&&
  -{1\over 2} \left( R_{i \overline{p}}S_{p \overline{j} k \overline{l}}  +
R_{p \overline{j}}S_{i \overline{p} k \overline{l}} + R_{k
\overline{p}}S_{i \overline{j} p \overline{l}} + R_{p
\overline{l}}S_{i \overline{j} k \overline{p}} \right)\\
&&   -{\mu \over 2} \left( R_{i \overline{p}} (g*g)_{p
\overline{j} k \overline{l}}  + R_{p \overline{j}} (g*g)_{i
\overline{p} k \overline{l}} + R_{k \overline{p}} (g*g)_{i
\overline{j} p \overline{l}} + R_{p \overline{l}}(g*g)_{i
\overline{j} k \overline{p}} \right)
\\ & = & \square S_{i \b j k\b l} +
 2 \mu S_{i \b j k\b l} + \mu
(S_{i\b j} g_{k \b l} + S_{k \b l} g_{ i \b j}) + \mu^2 (g*g)_{i\b
j p\b q} (g*g)_{q\b p k\b l}\\
&& \qquad  - 2 \mu S_{i \b j k\b l} - \mu (S_{i\b l k \b j} + S_{i
\b l k \b j }) - \mu^2 (g*g)_{i\b
p k \b q} (g*g)_{ p \b j q \b l}\\
&& \qquad + 2 \mu S_{i \b l k\b j} + \mu (S_{i\b l} g_{ k \b j} +
S_{k \b j} g_{i \b l}) + \mu^2 (g*g)_{i\b l p \b q} (g*g)_{ q \b p
k   \b j}\\
&& \qquad +  \mu (g*g)_{i \b j k \b l} -\mu \left( R_{i \b j} g_{k
\b l} + R_{k \b l} g_{i \b j} + R_{i \b l} g_{k \b j} + R_{k \b j}
g_{i \b l} \right)\\

& = & \square S_{i \b j k\b l} + \mu^2 \left((n+2) g_{i\b j}
g_{k\b l} + g_{i\b l } g_{k \b j} -2 g_{i\b j} g_{k\b l} - 2 g_{i
\b l} g_{k \b j} + (n+2) g_{i\b l} g_{k \b j} + g_{i \b j} g_{k \b
l} \right)\\
&& \qquad +  \mu (g*g)_{i \b j k \b l} +\mu \left( (S_{i\b j} -
R_{i \b j}) g_{k \b l} + (S_{k \b l} - R_{k \b l}) g_{i \b j} +
(S_{i\b l} - R_{i \b l}) g_{k \b j} + (S_{k \b j} - R_{k \b j})
g_{i \b l} \right)\\
& = & \square S_{i \b j k\b l}  + \mu \left((n+1)\mu + 1\right)
(g*g)_{i \b j k \b l}\\&& \qquad +\mu \left( (S_{i\b j} - R_{i \b
j}) g_{k \b l} + (S_{k \b l} - R_{k \b l}) g_{i \b j} + (S_{i\b l}
- R_{i \b l}) g_{k \b j} + (S_{k \b j} - R_{k \b j}) g_{i \b l}
\right)
\\ & = & \square S_{i \b j k\b l}  + \mu \left((n+1)\mu + 1\right)
(g*g)_{i \b j k \b l}\\&& \qquad -\mu (n+1) \mu \left( g_{i\b j}
g_{k \b l} + g_{k \b l}  g_{i \b j} + g_{i\b l} g_{k \b j} + g_{k
\b j} g_{i \b l} \right)\\
& = &  \square S_{i \b j k\b l}  + \mu \left(1 -(n+1)\mu \right)
(g*g)_{i \b j k \b l}
\end{array}
\]
Notice that $R_{i \b j  k \b l} = S_{i \b j k \b l} + \mu
(g*g)_{i\b j k\b l}$ and
\[\begin{array}{lcl} &&
{\p\over {\p t}} \left(\mu (g*g)_{i\b j k \b l}
\right)\\
 & = & \mu' (g*g)_{i \b j k \b l} + \mu \left(g_{i \b j}
(g_{k \b l} - R_{k \b l} ) + g_{i\b l} (g_{k \b j} - R_{k \b j}) +
g_{k \b l} (g_{i \b j} - R_{i \b j}) + g_{k \b j} (g_{i \b l} -
R_{i \b l})\right)\\
& = & \mu' (g*g)_{i \b j k \b l}  + 2 \mu (1 -(n+1)\mu) (g*g)_{i\b
j k \b l} - \mu (S_{i \b j} g_{k \b l} +g_{i \b j} S_{k \b l}+S_{i
\b l} g_{k \b j}+S_{k \b j} g_{i \b l})
\end{array}
\]
 we have
\[
\begin{array}{lcl} ({\p \over {\p t}} - \square) S_{i\b j k \b l}
& = & \left(\mu ((n+1)\mu - 1) - \mu'\right) (g*g)_{i \b j k \b
l}\\&& \qquad +\mu \left( S_{i\b j} g_{k \b l} + S_{k \b l} g_{i
\b j} + S_{i\b l} g_{k \b j} +  S_{k \b j} g_{i \b l} \right)\\
& = & \left(\mu ((n+1)\mu - 1) - \mu'\right) (g*g)_{i \b j k \b
l}\\&& \qquad +\mu \left( R_{i\b j} g_{k \b l} + R_{k \b l} g_{i
\b j} + R_{i\b l} g_{k \b j} +  R_{k \b j} g_{i \b l} \right)\\
&& - 2 \mu^2 (n+1) (g*g)_{i\b j k \b l}\\
& = & \left(\mu (-(n+1)\mu - 1) - \mu'\right) (g*g)_{i \b j k \b
l}\\&& \qquad +\mu \left( R_{i\b j} g_{k \b l} + R_{k \b l} g_{i
\b j} + R_{i\b l} g_{k \b j} +  R_{k \b j} g_{i \b l} \right).
\end{array} \]

 Note that we
assume that the positivity of the orthogonal bisectional curvature
is preserved by the K\"ahler Ricci flow.  Suppose $\mu(t) < 0$ is
the lower bound of the evolved bisectional curvature such that $S_{i
\b j k \b l} = R_{i \b j k \b l} - \mu(t) (g*g)_{i\b j k \b l} \geq
0. \;$ %We claim that $ -\mu(t)$ is a monotone decreasing function
%over any interval where the Ricci curvature is uniformly bounded.
%If the bound of Ricci curvature is {\it a priori} preserved over the
%entire flow, then the lower bound of the holomorphic sectional
%curvature will increase to $0$ as $t\rightarrow \infty.$
At the
point and time where $S_{i\b j k \b l}$ reaches $0,$ then this
minimum is realized by a holomorphic sectional curvature. Without
loss of generality, we may assume that it is $S_{1\b 1 1\b 1}(p,
t_0) = 0.\;$ Then at $t=t_0$,
\[
   R_{1\b 1 1\b 1} = 2\mu.
\]
We can diagonalize so that
\[
  R_{1 \b 1 k \b l}(p, t) = \lambda_k \delta_{k l}, \forall \;
  k,\;l = 1,2,\cdots n.
\]
Here $\lambda_1 =  2 \mu$ and \[ A = R_{1\b 1} (p,t) = \sum_{k =
1}^n\;\lambda_k = + 2\mu + \displaystyle \sum_{k=2}^n\; \lambda_k
\geq c.\]
 As in the Mok's argument \cite{Mok88},
we have
\[\begin{array}{lcl}
  {\p \over {\p t}} S_{1\b 1 1\b 1} \mid_{(p, t)}  & \geq &  S_{1\b 1k \b
  l} S_{\b k l \b 1 1} + \left(\mu (-(n+1)\mu - 1) - \mu'\right) (g*g)_{1 \b 1 1 \b 1} +   2 A \mu (g*g)_{1 \b 1 1 \b
 1}\\
 & = & \displaystyle \; \sum_{k=2}^n (\lambda_k - \mu)^2 + 2 \left(\mu (-(n+1)\mu - 1) -
 \mu'\right) + 4 A \mu\\
& \geq & {1\over {n-1}} \left(\displaystyle \sum_{k=2}^n \lambda_k -
(n-1) \mu\right)^2   + 2 \left(\mu (-(n+1)\mu - 1) -
 \mu'\right) + 4 A \mu \\ & = &
 {1\over {n-1}} (A -(n+1) \mu)^2  + 2 \left(\mu (-(n+1)\mu - 1) -
 \mu'\right) + 4 A \mu\\
 & = & {A^2 \over {n-1}} - 2\mu (1 - {{n-3}\over {n-1}}
A) + {{n+1}\over {n-1}} (3-n) \mu^2 - 2 \mu'.\end{array}
\]
From the last expression, it shows that for $n = 2, 3,$ the lower
bound of holomorphic sectional curvature is preserved (when the
absolute lower bound is big enough).  However, more should be true.
 In fact, when the minimum of $S_{i\b j k\b l}$ is achieved by
 holomorphic sectional curvature, by an argument of linear algebra
 (c.f.  Lemma 2.1  and inequality
\ref{technicalemma} below.  This appears explicitely in \cite{Cao92}
(which eventually is due to R. Hamilton)),  we can
 squiz a little more to obtain
\[\begin{array}{lcl}
  {\p \over {\p t}} S_{1\b 1 1\b 1} \mid_{(p, t)}  & \geq & 2 S_{1\b 1k \b
  l} S_{\b k l \b 1 1} + \left(\mu (-(n+1)\mu - 1) - \mu'\right) (g*g)_{1 \b 1 1 \b 1} +   2 A \mu (g*g)_{1 \b 1 1 \b
 1}\\
 & = & 2 \displaystyle \; \sum_{k=2}^n (\lambda_k - \mu)^2 + 2 \left(\mu (-(n+1)\mu - 1) -
 \mu'\right) + 4 A \mu\\
& \geq & {2 \over {n-1}} \left(\displaystyle \sum_{k=2}^n
\lambda_k - (n-1) \mu\right)^2   + 2 \left(\mu (-(n+1)\mu - 1) -
 \mu'\right) + 4 A \mu \\ & = &
 {2\over {n-1}} (A -(n+1) \mu)^2  + 2 \left(\mu (-(n+1)\mu - 1) -
 \mu'\right) + 4 A \mu\\
 & = & {2A^2 \over {n-1}} - 2\mu (1 + {4 \over {n-1}}
A) + {{4(n+1)}\over {n-1}}  \mu^2 - 2 \mu'\\
& = & {2\over {n-1}} (A - 2\mu)^2 - {8\mu^2 \over {n-1}} - 2\mu +
{{4(n+1)}\over {n-1}}  \mu^2 - 2 \mu' \\
& = & {2\over {n-1}} (A - 2\mu)^2 + 4\mu^2 - 2\mu - 2 \mu'
.\end{array}
\]
Consequently, the lower bound of holomorphic sectional curvature
is preserved and graduately improve to $0$ as $t\rightarrow \infty.$ No
assumption on Ricci curvature needed to prove the preservation of
lower bound here. However, if a positive lower bound of Ricci curvature
is preserved,  then the holomorphic sectional curvature shall
becomes positive after finite time.

\end{proof}
\begin{lem} Let $A, B, C$ are complex matrix and $A,C$ are hermitian
matrix.  Suppose that the real quadritic form
\[
Q = A_{i\b j} x^i x^{\b j}  + C_{i\b j} y^k y^{\b l} + 2 Re(B_{i j}
x^i y^j +  B_{i\b j} x^i y^{\b j}) \geq 0,\qquad \forall\; x, y \in
{\mathbb C}^n.
\]
Then, we have
\[
  A_{i \b j} C_{\b j i} \geq |B_{ij}|^2 + |B_{i\b j}|^2.
\]
\end{lem}
Now we indicate how we can apply this lemma to the preceding proof.
Suppose that $S_{i\b j k \b l}$ achieve minimum in the direction
${\p \over {\p z_1}} $ and ${\p \over {\p z_{\b 1}}}.\;$  Then, for
any
\[
x =\displaystyle \sum_{i=1}^n\; x^i {\p \over {\p z_i}}, \qquad y
=\displaystyle \sum_{i=1}^n\; y^i {\p \over {\p z_i}}
\]
we have
\[
  {\p^2 \over { \p \epsilon^2}} S({\p \over {\p z_1}} +  \epsilon x , {\p \over {\p z_{\b
1}}} +\epsilon \b x, {\p \over {\p z_1}} +  \epsilon y , {\p \over
{\p z_{\b 1}}} +\epsilon \b y) \mid_{\epsilon = 0} \geq 0.
\]
In other words, we have
\[
  S_{1\b 1 k \b l} y^k y^{\b l} + S_{i \b j 1 \b 1} x^i x^{\b j} + 2
Re( S_{1 \b k 1 \b j} x^{\b k} y^{\b j}  + S_{1 \b k j \b 1} x^{\b
k} y^j) \geq 0.
\]
Applying the lemma above, we have
\begin{eqnarray}
S_{1\b 1 k \b l} S_{\b k l \b 1 1} & \geq  S_{1\b k 1 \b j} S_{\b
1 k \b 1 j}  +  S_{1 \b k j \b 1} S_{\b 1 k \& b j 1}\\
& \geq  S_{1\b k 1 \b j} S_{\b 1 k \b 1 j}  +  S_{1 \b  1  j \b k}
S_{\b 1  1 \b j k}. \label{technicalemma}
\end{eqnarray}

 Now we give a proof to Theorem
 \ref{th:positivebisectionalcurvimprove}.

\begin{proof} First we assume that holomorphic sectional curvature is still negative somewhere in $M$ and
\[
   R_{i\b j } \geq \nu g_{i\b j} > 0.
\]
Following the notation of the previous proof in this section, we
have (at the minimum of the holomorphic sectional curvature)
\[\begin{array}{lcl}
  {\p \over {\p t}} S_{1\b 1 1\b 1} \mid_{(p, t)}  & \geq & 2 S_{1\b 1k \b
  l} S_{\b k l \b 1 1} + \left(\mu (-(n+1)\mu - 1) - \mu'\right) (g*g)_{1 \b 1 1 \b 1} +   2 A \mu (g*g)_{1 \b 1 1 \b
 1}\\
 & = & 2 \displaystyle \; \sum_{k=2}^n (\lambda_k - \mu)^2 + 2 \left(\mu (-(n+1)\mu - 1) -
 \mu'\right) + 4 A \mu\\
& \geq & {2 \over {n-1}} \left(\displaystyle \sum_{k=2}^n \lambda_k
- (n-1) \mu\right)^2   + 2 \left(\mu (-(n+1)\mu - 1) -
 \mu'\right) + 4 A \mu \\ & = &
 {2\over {n-1}} (A -(n+1) \mu)^2  + 2 \left(\mu (-(n+1)\mu - 1) -
 \mu'\right) + 4 A \mu\\
 & = & {2A^2 \over {n-1}} - 2\mu (1 + {4 \over {n-1}}
A) + {{4(n+1)}\over {n-1}}  \mu^2 - 2 \mu'\\
& = & {2\over {n-1}} (A - 2\mu)^2 - {8\mu^2 \over {n-1}} - 2\mu +
{{4(n+1)}\over {n-1}}  \mu^2 - 2 \mu' \\
& = & {2\over {n-1}} (A - 2\mu)^2 + 4\mu^2 - 2\mu - 2 \mu'\\
& \geq & {2\over {n-1}} \nu^2 - 2 \mu' .\end{array}
\]
The last inequality hold since $\mu < 0.\;$  This shows that at most
finite time, the holomorphic sectional curvature will become
positive everywhere.  \\

 Nowe we return to the case where the bisectional curvature is already positive.  By assumption, we have
\[
   R_{i\b j } \geq \nu g_{i\b j}, {1\over 2} < \nu < 1,
\]
and the bisectional curvature \[ R_{i\b j k \b l} \geq \mu(t)
(g*g)_{i\b j k \b l}, \qquad {\rm where }\; \mu(0) \geq 0.
\]
Then,
\[
\begin{array}{lcl} &&
{\p\over {\p t}} \left(\mu (g*g)_{i\b j k \b l}
\right)\\
 & = & \mu' (g*g)_{i \b j k \b l} + \mu \left(g_{i \b j}
(g_{k \b l} - R_{k \b l} ) + g_{i\b l} (g_{k \b j} - R_{k \b j}) +
g_{k \b l} (g_{i \b j} - R_{i \b j}) + g_{k \b j} (g_{i \b l} -
R_{i \b l})\right)\\
& \leq & \mu' (g*g)_{i \b j k \b l}  + 2 \mu (1 -v) (g*g)_{i\b j k
\b l}.
\end{array}
\]
Consequently, we have
\[
\begin{array}{lcl} ({\p \over {\p t}} - \square) S_{i\b j k \b l}
& \geq & \left(\mu (2 v- (n+1)\mu - 1) - \mu'\right) (g*g)_{i \b j
k \b l}.
\end{array} \]

Let \[ \mu(t) = {{C e^{a t}}\over {C e^{a t} + 1}} \cdot a, \qquad
{\rm where} \;  C = {{\mu(0)}\over {1 - \mu(0)}} \;{\rm and}\; a = {{2\nu -1}\over {n+1}}.
\]
Then
\[
({\p \over {\p t}} - \square) S_{i\b j k \b l} \geq 0.
\]
\end{proof}

\section{Proof of Theorem \ref{th:mainwithP}}

\begin{proof}
By Theorem \ref{th:riccibecomepositive}, Ricci curvature is
uniformly bounded from below over the entire flow.   According to
Perelman, the scalar curvature is uniformly bounded over the entire
flow.  Consequently, the Ricci curvature is uniformly bounded from
below and above.   Since orthogonal bisectional curvature is
positive and the holomorphic sectional curvature bounded from below
(Theorem \ref{th:holosectimprov}), then all bisectional curvature is
uniformly bounded over the entire K\"ahler Ricci flow. By Perelman's
non-local collapsing lemma when curvature is bounded, we imply that
the diameter is uniformly bounded from above.  Consequently, the
injectivity radius must have a uniform positive lower bound over the
flow. Suppose that $J$ is the underlying complex structure of $M.\;$
For any time sequence $\{t_i, i \in \NN\},$ there exists a
subsequence $\{g(t_i), i \in \NN\}$ and a sequence of diffeomorphism
$\rho_i: M \rightarrow M$ such that the sequence of pointed manifold
$\{(M, \rho_i^* g(t_i), J_i = \rho_i^*J), i\in \NN\}$ converges to a
K\"ahler Ricci soliton $g_\infty$ with non-negative bisectional
curvature in $M,$ perhaps with a different complex structure
$J_\infty$. Using Theorem \ref{th:holosectimprov}, the bisectional
curvature of the limit metric $g_\infty$ must be non-negative. Since
$g_\infty$ is a K\"ahler Ricci soliton, then either the Ricci
curvature of $g_\infty$ has no null direction at all or the K\"ahler
manifold will be split %holomorphically
 into a product of at least
two K\"ahler sub-manifolds. %Since the complex structure $J_\infty$
%of the limiting metric is adjacent to the origininal  complex
%structure $J$ (i.e., $\displaystyle \lim_{i \rightarrow \infty}\;
%J_i =J _\infty$),
Then,  $M$ must be %holomorphically
 decomposed as a
product of two K\"ahler Submanifolds. % for $J_i$ when $i$ is large
%enough. On the other hand, $J_i$ and $J$ are complete equivalent.
Since we assume that $M$ is irreducible, this leads a contradiction.
Thus, the Ricci curvature $R_{i\b j}(g_\infty)
>  0$ for any K\"ahler Ricci Soliton arised this way.  Since $M$ is compact, then there exists a small constant
$\epsilon_0 > 0$ such that $R_{i\b j}(g_\infty) > \epsilon_0.\;$ It
is then straightforward to see that there is a positive constant
$\epsilon_0> 0$ which doesn't depend on which time subsequences we
choose. In other words, for $t\geq t_0$ big enough, the Ricci
curvature is already positive; moreover, a positive uniform lower
bound is a priori preserved over the K\"ahler Ricci flow after $t
> t_0$. Following Theorem 1.5, the holomorphic sectional curvature must
become positive after finite time beyond $t= t_0.\;$   Since the
orthogonal bisectional curvature is always positive before
$t=\infty, $ we show that the evolved K\"ahler metric $g(t)$, after
at most finite time, must have positive bisectional curvature.
Appealing to the Frankel conjecture, then $M$ is
$\mathbb{C}\mathbb{P}^n.\;$ Then, the exponential convergence of the
flow
follows our earlier work \cite{chentian001} \cite{chentian002}.\\
\end{proof}
\section{Proof of Theorem \ref{th:mainwithP2}}
First, let's state a parabolic version of Moser iteration lemma as
in \cite{chen-12-04} without proof.

\begin{lem} \label{lem:ricciiteration} If the Poincare constant
and the Sobolev constant of the evolving K\"ahler metrics $g(t)$ are
both uniformly, and if a non-negative function $u$ satisfying the
following inequality
\[
{\p \over {\p t}} u \leq \triangle u + f(x, t)\; u, \qquad \forall
a< t < b,
\]
where  $|f|_{L^p(M,g(t))}$ is uniformly bounded by some constant
$c,$ then for any $\lambda \in (0,1)$ fixed, we have
\[\displaystyle \max_{(1-\lambda) a + \lambda b \leq t \leq b} u \leq C(c, b-a, \lambda)\int_a^b\; \int_M u.
  \]
\end{lem}

Now we prove Theorem \ref{th:mainwithP2}.
\begin{proof}

We want to prove that each of the three conditions will implies the
flow converges to a Fubni-study metrics exponentially. Naturally,
this break into three cases.\\

{\bf Case 1}.   Let us first suppose that $E_1$ has a lower bound.
Theorem \ref{th:riccibecomepositive} implies that $Ric (g(t))> -1$
hold after finite time.   The energy functional $E_1$ is monotonely
decreases afterwards.  Since $E_1$ is assumed to have a lower bound,
then there exists a $t_0 \geq 0$ such that

\[
\int_{t_0}^\infty \;d\, t \displaystyle\int_M\; |Ric(g(t)) -
\omega(g(t))|^2\; \omega(g(t))^n \leq \epsilon(n).
\]
 According to Theorem \ref{th:holosectimprov} and the fact that the
scalar curvature is uniformly bounded, the Riemannian curvature is
uniformly bounded from below and above on the entire K\"ahler Ricci
flow.  Therefore, the diameter are uniformly bounded on the K\"ahler
Ricci flow as in \cite{chentian002}.  Consequently, the injectivity
radius must have a uniformly positive lower bounded over the entire
flow. In other words, both Sobleve constant and Poincare constants
of the evolving metrics are uniformly bounded from above. According
to the iteration Lemma \ref{lem:ricciiteration}, the Ricci curvature
is uniformly pinched  towards identity by the constant $\epsilon(t)
\rightarrow 0$ such that
 we have
\[
   1-\epsilon(t)< Ric(g(t)) \leq 1+\epsilon(t), \qquad \forall t >
   t_0.
\]
In particular, let $\nu(t)$ to be the lower bound of $Ric(g(t)).\;$
Then $\displaystyle\lim_{t\rightarrow \infty} \nu(t) =  1.\;$
Appealing to Theorem \ref{th:positivebisectionalcurvimprove}, the
lower bound of the bisectional curvature eventually improve to
${1\over {n+1}}.\;$ However, at any point $p\in M$, adopting an
orthonomal frame, we have
\[\begin{array}{lcl}
  R_{i \b i} (p) & =&  \displaystyle \sum_{k = 1}^n \; R_{i\b i k \b
k}\\ & \geq &  \displaystyle \sum_{1\leq k \leq n; k \neq i} \;
R_{i\b i k \b k} + R_{i \b i i \b i}.
\end{array}
\]
Thus, if the lower bound of the bisectional curvature improve to
${1\over {n+1}}$ and the Ricci tensor approaches to the identity,
then the full bisectional curvature approach to the bisectional
curvature of a Fubini Study metric. Consequently, the underlying
manifold is $\mathbb{C}\mathbb{P}^n.\;$ Following Theorem
\cite{chentian002}, the flow converges exponentially fast
to a metric with constant bisectional curvature. \\

{\bf Case 2}.  Let us assume now that the Mabuchi energy has a
uniform lower bound.  Then, we have
\[\int_0^\infty\;d\,t\;\displaystyle \int_M\; |\nabla {{\p \varphi}\over {\p t}}
\mid_{\varphi(t)}^2 \; \omega_{\varphi(t)}^n \leq C.
\]
The evolution equation for $F(t) = |\nabla {{\p \varphi}\over {\p
t}} \mid_{\varphi(t)}^2$ is simple
\[
  {\p \over {\p t}} F = \triangle_{\varphi(t)} F + F - \left({{\p \varphi}\over {\p t}}\right)_{\alpha \b \beta} \cdot
  \left({{\p \varphi}\over {\p t}}\right)_{\b \alpha \beta} -  \left({{\p \varphi}\over {\p t}}\right)_{\alpha  \beta} \cdot
  \left({{\p \varphi}\over {\p t}}\right)_{\b \alpha \b \beta}.
\]
Following the iteration Lemma \ref{lem:ricciiteration}, we have
\[
 F = |\nabla {{\p \varphi(t)} \over {\p t}}|^2_{\varphi(t)} \rightarrow 0.
\]
Consider the evolution of $\int_M \; F $ over time.
\[\begin{array}{lcl}
{d\over {d\, t}} \displaystyle \int_M\; F\;\omega_{\varphi(t)}^n & =
&  \displaystyle \int_M\; F\;\omega_{\varphi(t)}^n + \displaystyle
\int_M \; F (n- R)
  \omega_{\varphi(t)}^n \\ && \qquad -
\displaystyle \int_M\; \left({{\p \varphi}\over {\p
t}}\right)_{\alpha \b \beta} \cdot
  \left({{\p \varphi}\over {\p t}}\right)_{\b \alpha \beta} \omega_{\varphi(t)}^n-   \displaystyle \int_M\; \left({{\p \varphi}\over {\p t}}\right)_{\alpha  \beta} \cdot
  \left({{\p \varphi}\over {\p t}}\right)_{\b \alpha \b \beta}
  \omega_{\varphi(t)}^n.
\end{array}
\]

Since the scalar curvature is uniformly bounded and $F\rightarrow 0$
pointwisely, we have

\[\begin{array}{l} \displaystyle \lim_{A\rightarrow \infty}
\int_A^{A+1}\;d\,t \;\displaystyle\int_M\; |Ric(g(t)) -
\omega(g(t))|^2\; \omega(g(t))^n  \\ \qquad \qquad  =  \displaystyle \lim_{A\rightarrow \infty}
\int_A^{A+1}\;d\,t \; \displaystyle \int_M\; \left({{\p
\varphi}\over {\p t}}\right)_{\alpha \b \beta} \cdot
  \left({{\p \varphi}\over {\p t}}\right)_{\b \alpha \beta} \omega_{\varphi(t)}^n
  \\  \qquad \qquad =   0.
  \end{array}
\]
Following the iteration Lemma again, we obtain
\[
   1-\epsilon(t)< Ric(g(t)) \leq 1+\epsilon(t), \qquad \forall t >
   t_0.
\]
Appealing to Theorem \ref{th:positivebisectionalcurvimprove}, the
bisectional curvature eventually pinch towards that of Fubini Study
metric.  Therefore the underlying manifold is
$\mathbb{C}\mathbb{P}^n.\;$ Following Theorem \cite{chentian002},
the flow converges exponentially fast to a metric with constant
bisectional curvature.\\

{\bf Case 3}.   Let us assume that  there exists a K\"ahler Einstein
metric in the canonical K\"ahler class.  By a well known theorem,
the Mabuchi energy has a uniform lower bound in this class. We can
reduce this to {\bf Case 2} before and  prove the exponential
convergence
of the K\"ahler Ricci flow.\\

We then complete our proof in all three cases.
\end{proof}
\section{Future problems}
The following questions are interesting:
\begin{enumerate}
\item Theorem \ref{th:mainwithP} is an extension of Frankel
conjecture.  Does there exists a direct proof as in Siu-Yau,
Morri's well known theorem?
 \item In this paper, we numerate quite a few new Maximum
 principle theorem. Is the lower bound of Ricci
 curvature is always preserved in K\"ahler Ricci flow? Perleman's work seems to suggest that this
 is true.
 \item It is known that the scalar curvature
 is bounded from above by exponential function over the K\"ahler Ricci flow.  Does the same hold
 for bisectional curvature?
\item How to derive the bound of the scalar curvature without
using Perelman's deep result? \item In this paper, we use
Perelman's result to derive scalar curvature bound.  We obtain
positive lower bound control on Ricci curvature after finite time
via an argument of taking limit.   Is that  possible to have a
more direct argument? \item Is compact K\"ahler manifold with
positive orthogonal bisectional curvature necessary have positive
first Chern class (except dimension $1$)? The answer should be
yes, at least when dimension is high enough. \item What happen in
the case of positive 2 curvature operator in the sense of H. Chen
\cite{[4]}.
\item Is there a (negative) lower bound of the holomophic sectional curvature which is preserved  under the K\"ahler Ricci flow  ?  This lower bound shall depends on the initial K\"ahler metric.  
My inclination  is that the answer shall be affirmative.
\end{enumerate}

%\bibliography{testa_bryant}
Xiuxiong Chen, Department of Mathematics/ University of Wisconsin/
Madison, WI 53706/ USA/ xiu@math.wisc.edu
\end{document}